\documentclass[1p]{elsarticle}
\usepackage{amssymb}
\usepackage{amsfonts}

\usepackage[polish, english]{babel}
\usepackage{polski}
\usepackage{tikz}

\newtheorem{theorem}{Theorem}
\newtheorem{conjecture}[theorem]{Conjecture}
\newtheorem{corollary}[theorem]{Corollary}
\newtheorem{lemma}[theorem]{Lemma}
\newtheorem{observation}[theorem]{Observation}

\newproof{pf}{Proof}

\begin{document}

\title{On asymptotically tight bound for the conflict-free chromatic index of nearly regular graphs}

\author[agh]{Mateusz Kamyczura}
\ead{kamyczuram@gmail.com}

\author[agh]{Jakub Przyby{\l}o}
\ead{jakubprz@agh.edu.pl}

\address[agh]{AGH University of Krakow, al. A. Mickiewicza 30, 30-059 Krakow, Poland}

\begin{abstract}
Let $G$ be a graph of maximum degree $\Delta$ which does not contain isolated vertices. An edge coloring $c$ of $G$ is called \emph{conflict-free} if each edge's closed neighborhood includes a uniquely colored element. The least number of colors admitting such $c$ is called the \emph{conflict-free chromatic index} of $G$ and denoted $\chi'_{\rm CF}(G)$. It is known that in general $\chi'_{\rm CF}(G)\leq 3 \lceil \log_2\Delta \rceil+1$, while there is a family of graphs, e.g. the complete graphs, for which $\chi'_{\rm CF}(G)\geq (1-o(1))\log_2\Delta$. In the present paper we provide the asymptotically tight upper bound  $\chi'_{\rm CF}(G)\leq (1+o(1))\log_2\Delta$  for regular and nearly regular graphs, which in particular implies that the same 
bound holds a.a.s. for a random graph $G=G(n,p)$ whenever $p\gg n^{-\varepsilon}$ for any fixed constant $\varepsilon\in (0,1)$. Our proof is probabilistic and exploits classic results of Hall and Berge. This was inspired by our approach utilized in the particular case of complete graphs, for which we give a more specific upper bound. We also observe that almost the same bounds hold in the open neighborhood regime.
\end{abstract}

\maketitle

\section{Introduction and notation}

Let $G=(V,E)$ be a graph with maximum degree $\Delta$ and minimum degree $\delta>0$.
By an \emph{edge coloring} $c$ of $G$ we understand any
assignment of colors to the edges of $G$.
Note we do not require $c$ to be proper. 
We denote by $E_G(v)$ the set of edges incident in $G$ with a vertex $v$, and we set
$E_G[uv]:=E_G(u)\cup E_G(v)$ to be the \emph{closed neighborhood} of any given edge $uv\in E$.
An edge $uv$ is called \emph{satisfied} by $c$ if there is at least one \emph{unique color} in its closed neighborhood, 
that is a color which appears exactly once in $E_G[uv]$. We give an advance notice to the fact we will be using the same notation when $c$ is a partial edge coloring (i.e. when $c$ attributes colors only to some subset of $E$). The coloring $c$ of $E$ is termed \emph{conflict-free} if all edges of $G$ are satisfied.
The minimum number of colors required to provide such a coloring is known as the \emph{conflict-free chromatic index} of $G$, and denoted $\chi'_{\rm CF}(G)$.

This concept was inspired by an earlier natural variant defined in the environment of vertex colorings.
More specifically, the least number of colors in a vertex coloring of $G$ assuring a unique color in the closed neighborhood of every vertex of $G$ is called the \emph{conflict-free chromatic number} of $G$ and denoted  $\chi_{\rm CF}(G)$. This graph invariant in turn stemmed from its correspondent considered in a more general setting of hypergraphs, applicable in a practical problem of 
non-interfering channel assignment in wireless networks, see
\cite{EvenEtAl,KostochkaEtAl,SmorodinskyApplications,SmorodinskyPhd}
for details and a few vital results.
In reference to the conflict-free chromatic number itself,
in 2009 Pach and Tardos~\cite{PachTardos} showed in particular that $\chi_{\rm CF}(G) = O(\ln^{2+\epsilon}\Delta)$. Later, the asymptotic behavior of this graph invariant (in the worst case) in terms of $\Delta$ was settled in two papers. Namely,  Bhyravarapu, Kalyanasundaram, and Mathew~\cite{Hindusi} improved the above upper bound to $\chi_{\rm CF}(G) = O(\ln^2\Delta)$, while Glebov, Szab\'o, and Tardos~\cite{GlebovEtAl} constructed a graph family with $\chi_{\rm CF}(G) = \Omega(\ln^2\Delta)$. 
Other directly related results, also concerning the open neighborhood setting, can additionally be found in~\cite{Hindusi3,Hindusi2,KellerEtAl,KostochkaEtAl}.
See also~\cite{CPS,EUN,EUN2,DWC-CHL,Fab,RH,MPK,CHL} for a list of related results within this intensively studied field, concerning similar concepts in the environment of proper colorings.

Note that investigating the conflict-free chromatic index $\chi'_{\rm CF}(G)$ may be regarded as studying $\chi_{\rm CF}(G)$, discussed above, within the family of line graphs.
The former graph invariant behaves however differently than the latter one. 
In particular, an upper bound of the form $\chi'_{\rm CF}(G)\leq C_1\log_2\Delta+C_2$, of a smaller order of magnitude than attainable in the case of $\chi_{\rm CF}(G)$, was provided in~\cite{DebskiPrzybylo}, where $C_1=9/\log_2\frac{1}{1-e^{-4}}\approx 337.5$ and $C_2$ are some universal constants.
The proof of this fact exploited the probabilistic method.
It was next improved by means of a simpler straightforward approach to the following best thus far known upper bound in terms of the maximum degree $\Delta$.
\begin{theorem}[\cite{MPK}] \label{MPK-Th}
    For every graph $ G$ without isolated vertices, $\chi'_{\rm CF}(G)\leq \lceil 3\log_2\Delta \rceil+1$.
\end{theorem}
This was in fact a consequence of a stronger upper bound, expressed in terms of the chromatic number $\chi(G)$ rather than the maximum degree, namely $\chi'_{\rm CF}(G)\leq \lceil 3\log_2\chi(G) \rceil+1$, cf.~\cite{MPK}. In~\cite{DebskiPrzybylo} D\k{e}bski and Przyby{\l}o also showed that in general this bound cannot be pushed much further down, demonstrating that one needs almost $\log_2\Delta$ colors in the case of the complete graphs, i.e. that $\chi'_{\rm CF}(G)\geq (1-o(1))\log_2\Delta$ for this graph family. It thus remained to potentially optimize the multiplicative constant in the upper bound from Theorem~\ref{MPK-Th}, which cannot be smaller than $1$.

We in particular prove that $1$ is in fact the right, and thus also optimal, multiplicative constant of the leading term $\log_2\Delta$ in the best possible upper bound for  $\chi'_{\rm CF}(G)$ in the case of regular graphs. 
Thus one of our main results can be formulated as follows.
\begin{theorem} \label{MPK-reg}
    For every nontrivial regular graph $ G$, $\chi'_{\rm CF}(G)\leq (1+o(1)) \log_2\Delta$.
\end{theorem}
In fact we will prove the assertion from Theorem~\ref{MPK-reg} above for a much wider family of nearly regular graphs, i.e. graphs with minimum degree
$\delta \geq \Delta - 2\sqrt{\Delta}\left(\ln\Delta\right)^{3/4}$, see Theorem~\ref{MainComplicated} below for details. Theorem~\ref{MPK-reg} is thus a direct corollary of the latter one. 
Another consequence of Theorem~\ref{MainComplicated}, proven for a broader family than just regular graphs, is the result concerning random graphs $G(n,p)$, for which we also prove that  $\chi'_{\rm CF}(G)\leq (1+o(1))\log_2\Delta$ a.a.s. for a wide spectrum of $p\gg n^{-\varepsilon}$, see Corollary~\ref{gnp-corollary} below.

The proof of our main result is based on several applications of the probabilistic method.
The crucial one of these exploits additionally the classic results of 
Hall and Berge on maximum matchings. 
This approach was inspired by our ideas designed to push down the best known upper bound $\chi_{\rm CF}'(K_n) \leq 2\lceil \log_2 n \rceil + 1$ from \cite{DebskiPrzybylo}  for the particular case of the complete graphs, which we improve to nearly optimal (or indeed optimal): $\chi_{\rm CF}'(K_n) \leq \lceil \log_2 (n-1) \rceil + 1$.

For convenience of a reader, in the next section we list the well-known theorems and tools
used within arguments throughout the paper. In the third section we 
demonstrate some of our ideas on the specific family of the complete graphs.
The following short section includes the main result accompanied by a couple remarks on on its proof,
which is then included in a separate section.
We close the paper with a concise concluding section,
including several conjectures and remarks concerning the open neighborhood setting.

\section{Tools}

We will use the basic symmetric variant of the Lov\'{a}sz Local Lemma, which can be found e.g. in \cite{LLL}. 

\begin{lemma}[Lov\'asz Local Lemma]\label{LLL}
Let $\Omega$ be a finite family of events in any probability space. Suppose that every event $A \in \Omega $ is mutually independent of a set of all the other events in $\Omega$ but at most $D$, and that $\mathbf{P}(A) \leq p$ for each $A \in \Omega$. If 
$$ep(D + 1) \leq 1,$$
then $\mathbf{P}(\bigcap_{A \in \Omega} \overline{A}) > 0$. 
\end{lemma}

We will moreover need two variants of the Chernoff Bounds.
The first one below can be found e.g. in~\cite{RandOm}, the second --  e.g. in~\cite{ChernoffBook}. 
\begin{lemma}[Chernoff Bound I]\label{Ch} 
Let $X= \sum_{i=1}^n X_i$ be a sum of independent random variables where $X_i=1$ with probability $p_i$, and $X_i = 0$ with probability $1 - p_i$. 
Then:
$$\mathbf{P}\left(X \leq \mathbf{E}(X)-t\right) \leq \exp\left(-\frac{t^2}{2\mathbf{E}(X)}\right)~~ {\rm for}  ~~ 0<t \leq \mathbf{E}(X).$$
\end{lemma}

\begin{lemma}[Chernoff Bound II]\label{Ch2} 
Let $X= \sum_{i=1}^n X_i$ be a sum of independent random variables where $X_i=1$ with probability $p_i$, and $X_i = 0$ with probability $1 - p_i$. Then for every $t>0$: 
      $$  \mathbf{P}\left( X \geq \mathbf{E}(X)+t\right) \leq \exp\left(-\frac{t^2}{t+2\mathbf{E}(X)}\right), $$ 
      hence
$ \mathbf{P}(X \geq \mathbf{E}(X)+t)\leq \exp(-\frac{t^2}{3\mathbf{E}(X)})$ if  $0<t \leq \mathbf{E}(X)$. 
\end{lemma}
It is straightforward to notice that the Chernoff Bounds above 
may also be applied for appropriate estimates on $\mathbf{E}(X)$. 
If $\mathbf{E}(X)\geq a>0$, then the assertion of Lemma~\ref{Ch} holds with $\mathbf{E}(X)$ replaced by $a$,
and similarly in the case of Lemma~\ref{Ch2}, when~$\mathbf{E}(X)\leq a$.

We will also make use of classic results: \emph{Halls's marriage theorem} and \emph{Berge's theorem on augmenting paths}.
Though these are very well known facts, we include them below for the sake of completeness,
specifying at the same time precisely the used notation. 
Recall that given a graph $G=(V,E)$ and a matching $M$ in it, we say a path $P$ of $G$ is \emph{alternating} if its consecutive edges alternate between $M$ and $E\smallsetminus M$.
We call it an \emph{augmenting path} (with respect to $M$) if it additionally starts and ends with unmatched vertices (i.e. not incident with $M$).
If $A,B\subseteq V$, we will denote by $G[A]$ the graph induced by $A$ in $G$, by $G[A,B]$ -- the graph induced by the edges joining $A$ and $B$ in $G$, and by $N_G(A)$, or simply $N(A)$, the set of neighbors of $A$, i.e. the vertices adjacent in $G$ with at least one vertex in $A$.
We also say that the matching $M$ \emph{saturates} $A$ if each vertex of $A$ is included in some edge of $M$.

\begin{theorem}[Hall's marriage theorem, \cite{Hall}]\label{HallTh}
Let $G=(V,E)$ be a bipartite graph with bipartition $V=X\cup Y$. Then
$G$ has a matching saturating $X$ if and only if 
\begin{equation}\label{HallsCond}
\forall S \subset X:~|N(S)|\geq |S|.
\end{equation}
\end{theorem}

\begin{theorem}[Berge's theorem, \cite{Berge}]\label{BergeTh}
A  matching $M$ in a graph $G$ is maximum
if and only if $G$ contains no augmenting path.
\end{theorem}

\section{Complete graphs}

In order to demonstrate some of our ideas we first prove that 
$\chi_{\rm CF}'(K_n)\leq \lceil \log_2(\Delta) \rceil +1$ 
for every nontrivial complete graph $K_n$.
\begin{theorem}\label{Th-ComplG}
For every $n\geq 2$,
$\chi_{\rm CF}'(K_n) \leq \lceil \log_2(n-1) \rceil+1$.
\end{theorem}

Note that to that end it is sufficient to provide a partial coloring $c$ of $K_n$ which satisfies all its edges and exploits at most $\lceil \log_2(n-1) \rceil$ colors (it is then sufficient to use one additional new color for the edges outside the domain of $c$). 
For this purpose we will in principle indicate a list of disjoint matchings in a given $K_n$, coloring each of these with one personal color. Each of these matchings will have roughly (at most) $n/4$ edges, and thus about $n/2$ vertices. Note every edge with one end in a vertex set of any given such matching $M$ and the other one outside it gets satisfied.  Thus one may expect any such matching to satisfy about half of the (remaining) edges incident with every vertex in our graph. After sequentially specifying roughly $\log_2n$ such matchings we may thus hope to satisfy all edges.
We use induction to prove that this is indeed the case. Theorem~\ref{Th-ComplG} is a direct corollary of the following slightly technical one,
within which we prove the existence of a desired partial coloring which additionally avoids the edges of any given matching in $K_n$ -- we call these edges \emph{blocked}.
\begin{theorem}\label{Th-ComplG-real}
For each complete graph $K_n$ with $n\geq 3$ and any blocked matching $M$ in it, there is a partial edge coloring of $K_n$ with at most $\lceil \log_2(n-1) \rceil$ colors which does not assign colors to the edges in $M$ and satisfies all edges of $K_n$. 
\end{theorem}

\begin{pf}
Let $K_n$ be a complete graph with vertex set $V$, $n\geq 3$, and let $M$ be a given blocked matching in $K_n$. 

Suppose first that $n\leq 6$.
Let $M'$ be a maximal matching in $K_n$ disjoint with $M$. Color the edges of $M'$ with pairwise distinct colors. As each edge of $K_n$ is adjacent with $M'$, all the edges are satisfied, 
and we have used $\lfloor\frac{n}{2}\rfloor \leq  \lceil \log_2(n-1) \rceil$ colors, as desired.

In order to inductively handle the remaining cases, let us first define auxiliary functions of a positive integer $n$, whose values are very close to $n/2$. Namely, let us set:
\begin{equation}\label{Hdefinition}
h(n)=2\left\lceil\frac{\left\lfloor\frac{n}{2}\right\rfloor}{2}\right\rceil;~~~~
\tilde{h}(n)=\max\left\{h(n),n-h(n)\right\}.
\end{equation}
Note that by definition, $h(n) \equiv 0\pmod 2$. 
One may easily verify that 
\begin{equation}\label{Hequality}
\tilde{h}(n) = \left\{\begin{array}{lll}
\left\lceil\frac{n}{2}\right\rceil & {\rm if} & n\not\equiv 2\pmod 4\\
\frac{n}{2}+1 & {\rm if} & n\equiv 2\pmod 4
\end{array} \right..
\end{equation}
Note that
\begin{equation}\label{HInequality}
1+\left\lceil\log_2\left(\tilde{h}(n)-1\right)\right\rceil \leq \left\lceil\log_2(n-1)\right\rceil 
\end{equation}
for every integer $n\geq 3$. Indeed, if $n\equiv 1\pmod 4$, 
then $\tilde{h}(n)=n-h(n)=\frac{n+1}{2}$, hence $2(\tilde{h}(n)-1)=n-1$ and thus~(\ref{HInequality}) follows. 
Otherwise, $\tilde{h}(n)=h(n)$, while $h(n)\leq \frac{n+1}{2}$ unless $n\equiv 2\pmod 4$. In the latter case however, $2(h(n)-1)=n$, but since $n-1$ cannot be an integer power of $2$ in such a case, (\ref{HInequality}) holds for all $n\geq 3$. 

Assume now that $n\geq 7$ and that our assertion holds for all complete graphs with orderers smaller than $n$.
Denote $M=\{a_1b_1, a_2b_2,\ldots,a_kb_k\}$, where $a_ib_i$ are pairwise distinct for all $i$. We may assume that $k=\lfloor\frac{n}{2}\rfloor$. Set $A=\{a_1,a_2,\ldots,a_k\}$, $B=\{b_1,b_2,\ldots,b_k\}$, $C=V\smallsetminus(A\cup B)$. (Note $|C|\leq 1$.)
Suppose we are able to partition $V$ to $V_1\cup V_2$ with $|V_1|=n_1$ and $|V_2|=n_2$ so that 
\begin{enumerate}
\item[(i)] $A\smallsetminus \{a_1\}\subseteq V_1$, $B\subseteq V_2$;
\item[(ii)] $n_1\equiv 0 \pmod 2$;
\item[(iii)] $\{n_1,n_2\}=\{h(n),n-h(n)\}$ and $n_1,n_2 \geq 3$.
\end{enumerate}
(Note that in order to assure (iii), 
it is sufficient to guarantee that $n_1=h(n)$ or $n_2=h(n)$.)
Then, by (i) no edge of $M$ is contained in $V_1$. By (ii) we may choose a matching $M''$ in $K_n[V_1]$ incident with all vertices in $V_1$ and color its edges with $1$. Note all edges in $K_n[V_1,V_2]$ become satisfied this way. By induction we may satisfy all edges of the complete graph $K_n[V_1]$
using at most $\lceil \log_2(n_1-1) \rceil$ new colors (different from $1$) to color some of its edges outside $M''$.
As by (i), $V_2$ includes at most one originally blocked edge (namely $a_1b_1$, if any), we may also satisfy all edges of the complete graph $K_n[V_2]$
using at most $\lceil \log_2(n_2-1) \rceil$  colors different from $1$ (note we may freely repeat colors used in $K_n[V_1]$) to color some of its edges distinct from $a_1b_1$. Consequently, all edges of $K_n$ get satisfied, while we have used at most $1+\max\{\lceil \log_2(n_1-1) \rceil, \lceil \log_2(n_2-1) \rceil\}$ colors, which, by (iii), (\ref{Hdefinition}) and(\ref{HInequality}), is at most  $\lceil \log_2(n-1) \rceil$.

It remains to prove that we may provide a partition of $V$ satisfying (i),(ii) and (iii). 
For  $n\not\equiv 2\pmod 4$ this is feasible as then $|A|=\lfloor\frac{n}{2}\rfloor\leq h(n)$ and $|A\cup C|\geq \lceil\frac{n}{2}\rceil \geq h(n)$, by~(\ref{Hequality}), and thus we may choose $V_1$ with $h(n)\equiv 0 \pmod 2$ elements so that $A\subseteq V_1\subseteq (A\cup C)$, and hence $B\subseteq V_2=V\smallsetminus V_1$.
For $n\equiv 2\pmod 4$ on the other hand, we may simply take $V_1=A\smallsetminus\{a_1\}$, as then $n_1=|V_1|=\frac{n}{2}-1\equiv 0 \pmod 2$ and $\{n_1,n_2\}=\{\frac{n}{2}-1,\frac{n}{2}+1\}=\{h(n),n-h(n)\}$, by~(\ref{Hequality}).
\qed
\end{pf}

\section{Main result}

We focus on nearly regular graphs of maximum degree $\Delta$ now. As it would be difficult to control some parameters within an inductive approach, and foremost very difficult to provide the base of induction, we modify our technique. 
This time we will generate about $\log_2\Delta$ matchings  randomly, all at once.
As in the case of the complete graphs, each of these matchings will have a personal color assigned and will span roughly $n/2$ vertices. Due to random choice, each of these will thus satisfy roughly half of all edges. In order to provide some independence between the matchings we will first partition the edges of our graph to 
roughly $\log_2\Delta$ subsets inducing near-regular subgraphs. Owing to this, each edge will have about $\log_2\Delta$ independent chances to get satisfied, and thus, with large probability, great majority of all edges incident with any given vertex of the graph will be satisfied after our randomized procedure, within which we exploit (color) only a small fraction of all edges. The remaining, unsatisfied edges will be handled with at the end by means of Theorem~\ref{MPK-Th}.
Let us also note that the required near-regularity of the graphs will be crucial for us while applying Hall's and Berge's theorems to generate desired matchings and within our randomized preparatory steps predating these. 

\begin{theorem}\label{MainComplicated}
There exists $\Delta_0$ such that  
$\chi_{\rm CF}'(G) \leq \log_2\Delta + 3\log_2\log_2\Delta +9$ for every graph $G$ with maximum degree $\Delta\geq \Delta_0$ and minimum degree $\delta \geq \Delta - 2\sqrt{\Delta}\left(\ln\Delta\right)^{\frac{3}{4}}$.
\end{theorem}

\section{Proof of Theorem~\ref{MainComplicated}}

\subsection{Preliminary graph decomposition}

Let $G=(V,E)$ be a graph with maximum degree $\Delta$ and minimum degree
\begin{equation}\label{delta-lower}
\delta \geq \Delta - 2\sqrt{\Delta}\left(\ln\Delta\right)^{\frac{3}{4}}
\end{equation}
where $\Delta$ is sufficiently large.  
Note we do not specify the lower bound $\Delta_0$
for $\Delta$, assuming throughout the proof it is large enough so that all explicite inequalities below involving $\Delta$ hold. We first decompose $G$ to 
$$s=\lceil \log_2\Delta \rceil$$ 
nearly regular subgraphs. To this end we independently assign every edge $e\in E$ to one of the sets 
$E_1,\ldots,E_s$ uniformly at random. We will show that with positive probability, for every resulting graph $G_i = (V,E_i)$ with minimum and maximum degrees $\delta_i$ and $\Delta_i$, respectively, 
\begin{equation}\label{ReqDeltaBounds}
\delta_i \geq \frac{\Delta}{s} - 3\sqrt{\Delta} ~~~~~ {\rm and} ~~~~~
      \Delta_i \leq \frac{\Delta}{s} + 3\sqrt{\Delta}.
      \end{equation}

Consider any vertex $v\in V$. 
Then for every $e\in E(v)=E_{G}(v)$ and each $i\in\{1,\ldots,s\}=[s]$, we denote
by $X_{i,v,e}$ the binary random variable equal to $1$ if $e \in E_i$ and $0$ otherwise. 	
Note $X_{i,v} = \sum_{e\in E(v)}X_{i,v,e}$ equals the degree $d_{G_i}(v)$ of $v$ in $G_i$.
Let us denote the following undesirable events.
\begin{itemize}
\item $U_{1,i,v}:~~X_{i,v} > \frac{\Delta}{s} + 3\sqrt{\Delta}$; 
\item $U_{2,i,v}:~~X_{i,v} < \frac{\Delta}{s} - 3\sqrt{\Delta}$. 
\end{itemize} 
Note that by~(\ref{delta-lower}), 
\begin{equation}\label{DegreExpBounds}
\frac{\Delta}{s}-\sqrt{\Delta} \leq \frac{\delta}{s} \leq \mathbf{E}(X_{i,v}) \leq  \frac{\Delta}{s}
\end{equation}
(for $\Delta$ large enough).
Since random variables $X_{i,v,e}$ are independent, we may apply the Chernoff Bounds, i.e. Lemmas \ref{Ch2} and \ref{Ch}, respectively, exploiting inequalities in~(\ref{DegreExpBounds}) and the fact that $\log_2\Delta > \ln\Delta$: 
\begin{equation} \label{U1bound}
  \mathbf{P}\left(U_{1,i,v}\right)
  \leq    \mathbf{P}\left(X_{i,v} \geq \frac{\Delta}{s} + 3\sqrt{\Delta}  \right) 
  \leq \exp\left(-\frac{9\Delta}{3\frac{\Delta}{s}}\right) \leq \Delta^{-3}
\end{equation}
\begin{equation} \label{U2bound}
  \mathbf{P}\left(U_{2,i,v}\right)
  \leq    \mathbf{P}\left(X_{i,v} \leq \left(\frac{\Delta}{s} - \sqrt{\Delta}\right) - 2\sqrt{\Delta}  \right) 
  \leq \exp\left(-\frac{4\Delta}{2\left(\frac{\Delta}{s}-\sqrt{\Delta}\right)}\right) \leq \Delta^{-2}.
\end{equation}

Each of the events $U_{1,i,v}, U_{2,i,v}$ is determined by random choices committed for edges incident with $v$, hence each of $U_{1,i,v}, U_{2,i,v}$ is mutually independent of all other events $U_{1,j,u}, U_{2,j,u}$ 
with $u$ at distance more than $1$ from $v$, i.e. of all other events but at most $D$ where
$D+1= 2\cdot s\cdot (\Delta+1)$.
Since, for $\Delta$ large enough,
$$e\cdot \Delta^{-2}\cdot \left(2\cdot s\cdot (\Delta+1)\right)<1,$$
by~(\ref{U1bound}), (\ref{U2bound}) and the Lov\'asz Local Lemma,
with positive probability none of the events $U_{1,i,v}, U_{2,i,v}$ with $i\in [s]$ and $v\in V$ occur.
There thus must exist an edge partition of $E$ and related graphs $G_i$ for which~(\ref{ReqDeltaBounds}) holds for every $i\in [s]$. 
We fix any such partition.

\subsection{Vertex decompositions}

Now in each $G_i$ we randomly choose a suitable matching. In fact, for every $i\in[s]$ we choose randomly only a certain partition of $V$ to two almost equal subsets with a few particular features. These, by Hall's and Berge's theorems will admit choosing a large matching $M$ in the first of the two subsets such that $M$ saturates almost all vertices from this subset. To guarantee this the mentioned first subset will be for technical reasons further partitioned to $3$ parts, say $V_1,V_2,V_3$, where $V_3$ will be very small. If one denotes the rest of the vertices by $V_4$, then the matching $M$, colored with its personal color will satisfy almost all edges between $V_1\cup V_2\cup V_3$ and $V_4$ (surely all between $V_1\cup V_2$ and $V_4$). Due to randomness, this guarantees that with high probability about half the edges incident with any given vertex is satisfied by this particular matching.
As our random choices for all $i\in [s]$ will be independent, great majority of all edges should get satisfied by $s$ disjoint matchings resulting from the randomized construction below.

For every $i\in [s]$, we denote:
\begin{equation}\label{EpsDef}
\varepsilon_i=\sqrt{\Delta_i}\ln\Delta_i. 
\end{equation}
Note that $\varepsilon_i\gg\sqrt{\Delta}$.
Let us associate with every vertex $v\in V$ exactly $s$ independent random variables $Z_{v,i}$, $i\in [s]$ such that 
\begin{equation}\label{Zdef}
Z_{v,i}=\left\{\begin{array}{lll}
1 & {\rm with~probability} & 1/4 - 3\varepsilon_i/\Delta_i\\ 
2 & {\rm with~probability} & 1/4\\ 
3 & {\rm with~probability} & 6 \varepsilon_i/\Delta_i\\ 
4 & {\rm with~probability} & 1/2 - 3 \varepsilon_i/\Delta_i 
\end{array}\right..
\end{equation}
These will be responsible for allocating vertices to subsets of designed vertex partitions in $G_i$'s and will help us to control edge distributions between these subsets, which are crucial for the sake of applications of Hall's and Berge's theorems. 
We precisely set:
$$V_q^{(i)}=\{v\in V~:~Z_{v,i}=q\}$$
for $i\in[s]$ and $q\in[4]$, which provides a vertex partition $V=V_1^{(i)}\cup V_2^{(i)}\cup V_3^{(i)}\cup V_4^{(i)}$ for every $G_i$. We also define random variables which for each $v\in V$ correspond to the number edges joining it with $V_1^{(i)}$, $V_2^{(i)}$ and $V_1^{(i)}\cup V_3^{(i)}$, respectively, in $G_i$:
$$Y_{v,1}^{(i)} = |V_1^{(i)}\cap N_{G_i}(v) |,~~~~
Y_{v,2}^{(i)} = |V_2^{(i)}\cap N_{G_i}(v)|,~~~~
Y_{v,13}^{(i)} = |(V_1^{(i)}\cup V_3^{(i)})\cap N_{G_i}(v)|.$$
Consider the following events for every $v\in V$ and each $i\in [s]$:
\begin{itemize}
    \item  $S_{v,1}^{(i)} :~~ Y^{(i)}_{v,1} < \frac{\Delta_i}{4} - 2\varepsilon_i $,
    \item  $B_{v,2}^{(i)} :~~ Y^{(i)}_{v,2} > \frac{\Delta_i}{4} - 2\varepsilon_i $,
    \item  $B_{v,13}^{(i)} :~~ Y^{(i)}_{v,13} > \frac{\Delta_i}{4} + \varepsilon_i $,
    \item  $S_{v,2}^{(i)} :~~ Y^{(i)}_{v,2} < \frac{\Delta_i}{4} + \varepsilon_i $.
\end{itemize}
We prove that the probability of each of these not to occur is very small. 

Note first that by~(\ref{Zdef}),
$\mathbf{E}(Y_{v,1}^{(i)})\leq\frac{\Delta_i}{4}-3\varepsilon_i$. Thus, by the Chernoff Bound, (\ref{EpsDef}) and~(\ref{ReqDeltaBounds}),
\begin{eqnarray}
\mathbf{P}\left(\overline{S_{v,1}^{(i)}}\right)
&=& \mathbf{P}\left(Y^{(i)}_{v,1} \geq \left(\frac{\Delta_i}{4} - 3\varepsilon_i\right)+\varepsilon_i\right)
\leq \exp\left(-\frac{\varepsilon_i^2}{3\left(\frac{\Delta_i}{4} - 3\varepsilon_i\right)}\right) \nonumber\\
&\leq& \exp\left(-\frac{4}{3}\ln^2\Delta_i\right)
 \leq \Delta^{-4}. \label{Sv1Bound}
\end{eqnarray}

Next, by~(\ref{Zdef}), (\ref{EpsDef}) and~(\ref{ReqDeltaBounds}), $\mathbf{E}(Y_{v,2}^{(i)})\geq\frac{\delta_i}{4} \geq\frac{\Delta_i}{4}-\varepsilon_i$. Thus, by the Chernoff Bound, (\ref{EpsDef}) and~(\ref{ReqDeltaBounds}),
\begin{eqnarray}
\mathbf{P}\left(\overline{B_{v,2}^{(i)}}\right)
&=& \mathbf{P}\left(Y^{(i)}_{v,2} \leq \left(\frac{\Delta_i}{4} - \varepsilon_i\right) - \varepsilon_i\right)
\leq \exp\left(-\frac{\varepsilon_i^2}{2\left(\frac{\Delta_i}{4} - \varepsilon_i\right)}\right) \nonumber\\
&\leq& \exp\left(-2\ln^2\Delta_i\right)
 \leq \Delta^{-4}. \label{Bv2Bound}
\end{eqnarray}

Subsequently, again by~(\ref{Zdef}), (\ref{EpsDef}) and~(\ref{ReqDeltaBounds}), 
$\mathbf{E}(Y_{v,13}^{(i)})\geq \delta_i(\frac{1}{4}+\frac{3\varepsilon_i}{\Delta_i}) 
\geq (\Delta_i-6\sqrt{\Delta})(\frac{1}{4}+\frac{3\varepsilon_i}{\Delta_i}) 
\geq\frac{\Delta_i}{4}+2\varepsilon_i$. Thus, by the Chernoff Bound, (\ref{EpsDef}) and~(\ref{ReqDeltaBounds}),
\begin{eqnarray}
\mathbf{P}\left(\overline{B_{v,13}^{(i)}}\right)
&=& \mathbf{P}\left(Y^{(i)}_{v,13} \leq \left(\frac{\Delta_i}{4} + 2\varepsilon_i\right)-\varepsilon_i\right)
\leq \exp\left(-\frac{\varepsilon_i^2}{2\left(\frac{\Delta_i}{4} + 2\varepsilon_i\right)}\right) \nonumber\\
&\leq& \exp\left(-\ln^2\Delta_i\right)
 \leq \Delta^{-4}. \label{Bv13Bound}
\end{eqnarray}

Analogously, by~(\ref{Zdef}),
$\mathbf{E}(Y_{v,2}^{(i)})\leq\frac{\Delta_i}{4}$. Thus, by the Chernoff Bound, (\ref{EpsDef}) and~(\ref{ReqDeltaBounds}),
\begin{eqnarray}
\mathbf{P}\left(\overline{S_{v,2}^{(i)}}\right)
&=& \mathbf{P}\left(Y^{(i)}_{v,2} \geq \frac{\Delta_i}{4} +\varepsilon_i\right)
\leq \exp\left(-\frac{\varepsilon_i^2}{3\frac{\Delta_i}{4}}\right) \nonumber\\
&=& \exp\left(-\frac{4}{3}\ln^2\Delta_i\right)
 \leq \Delta^{-4}. \label{Sv2Bound}
\end{eqnarray}

Later we will explain why $S_{v,1}^{(i)}$,
$B_{v,2}^{(i)}$,
$B_{v,13}^{(i)}$,
$S_{v,2}^{(i)}$ guarantee existence of $s$ large and convenient matchings. 
We however cannot know a priori for any given vertex $v$ to which subset of a random partition of any $G_i$ this vertex will belong, and in particular whether $v$ will be incident with the resulting matching in $G_i$ or not. 
Thus, in order to show that with large probability many edges incident with any vertex $v$ will eventually be satisfied, we will prepare ourselves for every possible scenario (each of which will be be encoded by means of a binary vector $t$ below). To that end we introduce one more type of technical events, preceded by a few preparatory estimations. 

Let $t\in\{0,1\}^s$ be any binary vector of length $s$. For any $v\in V$, $i\in[s]$ and a neighbor $u\in N_{G}(v)$ of $v$ in $G$, we first define the following event:
$$D_{v,u,t,i}:~~ \left(t[i]=0  \ \wedge \ u \in V_4^{(i)} \cup V_3^{(i)}\right) \vee \left(t[i]=1  \ \wedge \ u \in V_1^{(i)} \cup V_2^{(i)} \cup V_3^{(i)}\right).$$
Note that for the given $t$ and $i$, exactly one of the two must hold: either $t[i]=0$ or $t[i]=1$.
Thus, the probability of $D_{v,u,t,i}$ to occur equals the probability of an appropriate of the two events whose alternative defines $D_{v,u,t,i}$ to appear. However,  by~(\ref{Zdef}), the probability that $u$ belongs to $V_4^{(i)} \cup V_3^{(i)}$ is the same as the probability that it belongs to $V_1^{(i)} \cup V_2^{(i)} \cup V_3^{(i)}$, hence by~(\ref{Zdef}),  (\ref{EpsDef}) and~(\ref{ReqDeltaBounds}),
\begin{equation}\label{DvutiProb}
\mathbf{P}\left(D_{v,u,t,i}\right)=\frac{1}{2}+\frac{3\varepsilon_i}{\Delta_i}
= \frac{1}{2}+\frac{3\ln\Delta_i}{\sqrt{\Delta_i}} 
< \frac{1}{2}+\frac{s^2}{\sqrt{\Delta}}.
\end{equation}

Next, for the given $v$, $t$ and $u\in N_{G}(v)$, we set $D_{v,u,t}=\bigcap_{i\in[s]}D_{v,u,t,i}$, i.e.,
$$D_{v,u,t}:~~~~  \forall_{i\in [s]}: \left(\left(t[i]=0  \ \wedge \ u \in V_4^{(i)} \cup V_3^{(i)}\right) \vee \left(t[i]=1  \ \wedge \ u \in V_1^{(i)} \cup V_2^{(i)} \cup V_3^{(i)}\right)\right).$$
As choices for distinct $i\in [s]$ are independent, by~(\ref{DvutiProb}),
\begin{equation}\label{DvutProb}
\mathbf{P}\left(D_{v,u,t}\right) <  \left( \frac{1}{2} + \frac{s^2}{\sqrt{\Delta}}\right)^{s} 
           \leq \left(\frac{1}{2}\right)^s\left(1+\frac{1}{s}\right)^s 
           \leq \frac{e}{\Delta}.
\end{equation}
Let now $I_{v,u,t}$ be a binary random variable which takes value $1$ if $D_{v,u,t}$ holds and $0$ otherwise.
Set $I_{v,t}=\sum_{u\in N_G(v)}I_{v,u,t}$.
We will show that we may make our choices determining vertex partitions so that for every $v\in V$ and $t\in\{0,1\}^s$ the following event holds:
$$D_{v,t}:~~ I_{v,t} < 3\log_2\Delta,$$ 
i.e. that for less than $3\log_2\Delta$ 
neighbors $u$ of $v$ in $G$, $D_{v,u,t}$ holds.
Note that since choices for all $u\in N_G(v)$ are independent and by~(\ref{DvutProb}), we have $\mathbf{E}(I_{v,t})<e$, then by the Chernoff Bound from Lemma~\ref{Ch2}, since $\log_2\Delta > 1.4\ln\Delta$,
\begin{equation}\label{DvBound} 
\mathbf{P}\left(\overline{D_{v,t}}\right) = \mathbf{P}\left( I_{v,t} \geq e +\left(3\log_2\Delta-e\right)\right)
\leq \exp\left(-\frac{\left(3\log_2\Delta-e\right)^2}{\left(3\log_2\Delta-e\right)+2e}\right)
\leq \Delta^{-4.1}. 
\end{equation}

Finally we associate with every $v\in V$ an aggregated undesirable event that at least one of the events $S_{v,1}^{(i)}$,
$B_{v,2}^{(i)}$,
$B_{v,13}^{(i)}$,
$S_{v,2}^{(i)}$ or $D_{v,t}$ does not hold, i.e. we set:
$$A_v:~~\bigvee_{i\in[s]}\left(\overline{S_{v,1}^{(i)}}\vee
\overline{B_{v,2}^{(i)}}\vee
\overline{B_{v,13}^{(i)}}\vee
\overline{S_{v,2}^{(i)}}\right) \vee\bigvee_{t\in\{0,1\}^s}\overline{D_{v,t}}$$
Note that as there are $2^s<2\Delta$ distinct $t\in \{0,1\}^s$ and $s$ potential values of each $i$, by~(\ref{Sv1Bound}), (\ref{Bv2Bound}), (\ref{Bv13Bound}), (\ref{Sv2Bound}) and~(\ref{DvBound}), we have that 
\begin{equation}\label{AvBound}
\mathbf{P}\left(A_v\right) \leq 4s\Delta^{-4} + 2\Delta\cdot\Delta^{-4.1} \leq \Delta^{-3}. 
\end{equation}

Note that each $A_v$ is determined by  random choices associated with the neighbors of $v$.
Thus, every $A_v$ is  is mutually independent of all other events $A_u$ 
with $u$ at distance more than $2$ from $v$, i.e. of all other but at most $\Delta^2$ events.
By~(\ref{AvBound}) and the Lov\'asz Local Lemma,
with positive probability none of the events $A_v$ occur.
There thus must exist a set of partitions $V=V_1^{(i)}\cup V_2^{(i)}\cup V_3^{(i)}\cup V_4^{(i)}$ for all $i\in [s]$ such that $S_{v,1}^{(i)}$,
$B_{v,2}^{(i)}$,
$B_{v,13}^{(i)}$,
$S_{v,2}^{(i)}$ and $D_{v,t}$ hold for all $v\in V$, $i\in [s]$ and $t\in \{0,1\}^s$.
We fix any such collection of vertex partitions.

\subsection{Existence of suitable matchings}

We now argue that we may choose in every $G_i$ a matching $M_i$ with all edge ends in $V_1^{(i)}\cup V_2^{(i)}\cup V_3^{(i)}$ which saturates all vertices in $V_1^{(i)}\cup V_2^{(i)}$. Consider first a bipartite graph $H_i=G_i[V_1^{(i)}, V_2^{(i)}]$. Note that for all $u\in V_1^{(i)}$ and $v\in V_2^{(i)}$, by  $S_{u,1}^{(i)}$ and $B_{v,2}^{(i)}$, we have $d_{H_i}(u) > d_{H_i}(v)$.
It is straightforward to notice that  this guarantees via a simple double edge counting that Hall's condition~(\ref{HallsCond}) is satisfied, and thus by Theorem~\ref{HallTh} there is a matching $M'_i$ saturating $V_1^{(i)}$ in $H_i$.
Such $M'_i$ is obviously also included in the bipartite supergraph $H'_i=G_i[V_1^{(i)}\cup V_3^{(i)}, V_2^{(i)}]$ of $H_i$. 
This time, by  $S_{u,2}^{(i)}$ and $B_{v,13}^{(i)}$, we have $d_{H'_i}(u) > d_{H'_i}(v)$ for all $u\in V_2^{(i)}$ and $v\in V_1^{(i)}\cup V_3^{(i)}$.
Hence, analogously as above, Hall's theorem guarantees that every maximum matching in $H'_i$ saturates $V_2^{(i)}$. Thus, if this is not the case for $M'_i$, then $M'_i$ is not a maximum matching in $H'_i$, and hence, by Berge's theorem -- Theorem~\ref{BergeTh}, $H'_i$ has an augmenting path $P'_i$ containing $M'_i$. 
A symmetric difference of $P'_i$ and $M'_i$ yields a matching $M''_i$ which is larger than $M'_i$ and saturates the ends of all edges in $M'_i$, in particular all vertices in $V_1^{(i)}$. We may continue building such size increasing matchings along consecutive augmenting paths until we get a maximum matching $M_i$ in $H'_i$, which by construction must saturate $V_1^{(i)}$, while by Berge's theorem must also saturate $V_2^{(i)}$.
Thus, $M_i$ is indeed our desired matching.

For each $i\in[s]$ we fix any such matching $M_i$ and color all its edges with $i$.
We denote the resulting partial coloring of $G$ by $c$.
Consider any given vertex $v\in V$. Let $t_v$ be a vector in $\{0,1\}^s$ defined as follows:
$$t_v[i]=\left\{
\begin{array}{lll}
1 & {\rm if} & M_i~ {\rm is~incident~with}~ v\\
0 & {\rm if} & M_i~ {\rm is~not~incident~with}~ v
\end{array}\right..$$ 
As $D_{v,t}$ holds for each $t\in \{0,1\}^s$, hence in particular for $t=t_v$,
it follows that for all but at most $3\log_2\Delta$ vertices $u$ in $N_G(v)$ the event $D_{v,u,t_v}$ is not fulfilled. 
We note that for any such $u$ (for which $D_{v,u,t_v}$ does not hold) the edge $uv$ is already satisfied by our partial coloring $c$. Indeed, $\overline{D_{v,u,t_v}}$ implies that there is $i'\in [s]$ for which $D_{v,u,t_v,i'}$ is not true. If $t_v[i']=0$, this implies that $u\notin V_4^{(i')}\cup V_3^{(i')}$, i.e. $u\in V_1^{(i')}\cup V_2^{(i')}$, which means that $M_{i'}$ is incident with $u$. As $t_v[i']=0$ certifies that $M_{i'}$ is not incident with $v$, it follows that $uv$ is satisfied by $c$. On the other hand, if $t_v[i']=1$, i.e. $M_{i'}$ is incident with $v$, then $\overline{D_{v,u,t_v,i'}}$ implies that $u\notin V_1^{(i')}\cup V_2^{(i')}\cup V_3^{(i')}$, i.e. $u\in V_4^{(i')}$,
which means that $M_{i'}$ is not incident with $u$. Consequently, also in this case $uv$ is satisfied by $c$.

\subsection{Coloring the remaining edges}\label{SubSecFinalCol}

By the reasoning above, every vertex in $G$ is indent with at most 
$3\log_2\Delta$  yet not satisfied edges. Note these could not have been colored thus far.
Let us dente by $G'$ the graph induced by all these edges, hence 
$\Delta(G')\leq 3\log_2\Delta$. 
By Theorem~\ref{MPK-Th}, it is thus sufficient to use at most 
$\lceil 3\log_2(3\log_2\Delta) \rceil+1$  new colors in order to provide an edge coloring of $G'$
satisfying all its edges. 
We then color all yet uncolored edges in $G$ with one more final color, which completes the proof, as in total we have used no more than
$$s+\left(3\log_2(3\log_2\Delta)+2\right)+1 < \log_2\Delta + 3\log_2\log_2\Delta+9$$
colors.
\qed

\section{Random graphs}

Apart from directly implying Theorem~\ref{MPK-reg},
Theorem~\ref{MainComplicated} also applies within random graphs model via the following standard observation.
We include its proof for the sake of completeness. See e.g.~\cite{gnp} for all detailed definitions concerning this model and similar arguments.

\begin{observation}\label{gnp_obs}
    If $G= G(n,p)$ is a random graph with $p\gg n^{-\varepsilon}$ where $\varepsilon \in (0,1)$ is a constant, then a.a.s. 
    $\delta(G) \geq \Delta(G) - 2\sqrt{\Delta(G)}(\ln\Delta(G))^{3/4}$ and $\Delta(G) \gg n^{1-\varepsilon}$. 
\end{observation}

\begin{pf}
Given $p$ consistent with the assumption,
we denote by $d(v)$ the degree of any vertex $v$ in our random model. 
As ech potential edge $uv$ appears in $G(n,p)$ independently with probability $p$, we have that $\mathbf{E}(d(v))=(n-1)p\gg n^{1-\varepsilon}$. By the Chernoff Bound we thus obtain
\begin{eqnarray}
&&\mathbf{P}\left(\left|d(v)-\mathbf{E}(d(v))\right|\geq \sqrt{\mathbf{E}(d(v))}\left(\ln\mathbf{E}(d(v))\right)^\frac{3}{4}\right) 
\leq 2 \exp \left(-\frac{\left(\ln\mathbf{E}(d(v))\right)^\frac{3}{2}}{3}\right) \nonumber\\
&\leq& 2 \exp \left(-\frac{\left(1-\varepsilon\right)^\frac{3}{2}\left(\ln n\right)^\frac{3}{2}}{3}\right) 
= 2n^{-\frac{\left(1-\varepsilon\right)^\frac{3}{2}\sqrt{\ln n}}{3}} \nonumber
= o(n^{-1}).
\end{eqnarray}
Thus, the probability that the degree of at least one of $n$ vertices of our random graph deviates by more than $\sqrt{(n-1)p}(\ln((n-1)p))^{3/4}$ from $(n-1)p$ tends to $0$ with $n$. Hence, for $d=(n-1)p$ and the functions $f(x)=\sqrt{x}(\ln x)^{3/4}$ and $g(x)=x-2f(x)$, which are increasing for $x$ large enough,  a.a.s. the maximum degree $\Delta$ and the minimum degree $\delta$ of our graph fulfill:
\begin{eqnarray}
\Delta &\leq& d + f(d), \nonumber\\
\delta &\geq& d - f(d) = \left(d+f(d)\right)-2f(d) \nonumber\\
&\geq&  \left(d+f(d)\right)-2f(d+f(d))  
= g(d+f(d))
\geq g(\Delta), \nonumber
\end{eqnarray}
which completes the proof.
\qed
\end{pf}

Observation~\ref{gnp_obs} and Theorem~\ref{MainComplicated} imply the following corollary.
\begin{corollary}\label{gnp-corollary}
    If $G= G(n,p)$ is a random graph with $p\gg n^{-\varepsilon}$ where $\varepsilon \in (0,1)$ is a constant, then a.a.s.
    $\chi_{\rm CF}'(G) \leq \log_2\Delta + O(\log_2\log_2\Delta) = (1+o(1))\log_2\Delta$,
     where $\Delta$ denotes the maximum degree of $G$. 
\end{corollary}

\section{Concluding remarks}

We first refer to a naturally defined correspondent of $\chi_{\rm CF}'(G)$ where we require a uniquely colored element in the open neighborhood of every edge, rather than the closed one. We denote it by $\chi_{\rm OCF}'(G)$ for every graph $G$ without isolated vertices and without isolated edges. 
In case of the complete graphs one may mimik the proof of Theorem~\ref{Th-ComplG-real}, with a few modifications to obtain a similar upper bound as for $\chi_{\rm CF}'(G)$. Note we need to additionally guarantee that all edges of the matchings colored within our inductive approach are satisfied (as we are not certain of this within the open neighborhood regime).
This leads to a slightly more annoying small cases to analyze.
Consequently, some new exceptions emerge, in particular one, starting with 
$n=6$, which turns into an infinite family, regarding all orders of the form $n=2^q-2$, $q\geq 3$.
Also all orders $n=2^q+1$, $q\geq 3$, similarly as the ones above, seem to require one more color than in the case of the closed neighborhoods.
These two families additionally complicates our inductive approach slightly, as we need to strive to avoid certain orders while partitioning the vertex set to two subsets inducing complete subgraphs. The corresponding modifications are however achievable and straightforward.
We omit tedious details, stating only the result, which implies the same bound as in the case of  $\chi_{\rm CF}'(G)$ in most of the cases.
\begin{theorem}\label{Th-ComplG-open}
For every $n\geq 3$,
$\chi_{\rm OCF}'(K_n) \leq \lceil \log_2n \rceil+1$ if $n+2$ is not a natural power of $2$, and $\chi_{\rm OCF}'(K_n) \leq \lceil \log_2n \rceil+2$ otherwise.
\end{theorem}
We believe it would be interesting to finally settle the precise value of the both parameters in case of the complete graph. We suspect that the upper bounds in Theorems~\ref{Th-ComplG} and~\ref{Th-ComplG-open} are accurate or almost accurate.

As for our general bound from Theorem~\ref{MainComplicated}, it is enough to introduce some changes and supplementations in the very last part of the argument, only after providing our randomly constructed matchings $M_1,\ldots,M_s$. Note that almost all edges satisfied via these matchings remain also satisfied within the open neighborhood setting, except possible these in the matchings themselves. There can however be at most $s=\lceil\log_2\Delta\rceil$ of these incident with every vertex $v$ of the considered graph. Counting in the remaining at most $3\log_2\Delta$ potentially unsatisfied edges in $E(v)$, handled with in Subsection~\ref{SubSecFinalCol}, we are left with a new subgraph $G'$ of $G$ of unsatisfied edges of maximum degree $\Delta'\leq 6\ln\Delta$. There are many ways one may satisfy these limited number of edges incident with each vertex. One way is to use the approach exploited in~\cite{DebskiPrzybylo} to prove Theorem 1, implying that for any fixed $\alpha>0$ one may provide a \emph{vertex} coloring of any graph $H$ with $\delta(H)\geq \alpha \Delta(H)$ witnessing that $\chi_{\rm CF}(H)=O(\ln \Delta(H))$. The same can be proven 
in the open neighborhood regime even if additionally we are not allowed to color vertices from some set, say $A$, such that each vertex in $H$ has at most $\beta\Delta(H)$ neighbors in $A$, where
$\beta<\alpha$ is some fixed constant. The main, and essentially the sole, change one needs to introduce to the mentioned proof of Theorem 1 from~\cite{DebskiPrzybylo} is to avoid vertices in $A$ while randomly choosing vertices in the first part of this proof. In order to apply this result for our purposes, we may first randomly supplement $G'$ with some of the remaining uncolored edges from $G$ (including each of these independently with probability $p'=\Theta(\Delta^{-1}\ln\Delta)$) so that the resulting supergraph $G''$ of $G'$ has minimum degree $\delta(G'')=\Omega(\ln\Delta)$ and maximum degree $\Delta=O(\ln\Delta)$. The same will thus hold for the line graph $H=L(G'')$ of such $G''$, so we may apply the mentioned modification  of Theorem 1 from~\cite{DebskiPrzybylo} to this $H$, with $A$ comprising of elements of $M_1,\ldots,M_s$ from $G''$ (if we choose $p'$ large enough to guarantee that $\beta<\alpha$), using a new set of colors.
Consequently, it is sufficient to use $s+O(\ln\ln\Delta)$ colors to obtain a partial edge coloring of a given nearly regular graph  with a uniquely colored edge in the open neighborhood of every edge. We thus obtain the following, which also implies an obvious correspondent of Corollary~\ref{gnp-corollary} in the random model.
\begin{theorem}\label{Main-open}
Let $G$ be a graph of maximum degree $\Delta$ without isolated vertices and edges.
Then, 
$\chi_{\rm OCF}'(G) \leq \log_2\Delta + O(\log_2\log_2\Delta) = (1+o(1))\log_2\Delta$.
\end{theorem}

We conclude by posing several conjectures, two of which are just weaker variants of the other two.
\begin{conjecture}\label{ConjectureStrong-closed}
There is a constant $C$ such that for every graph $G$ of maximum degree $\Delta$ without isolated vertices,
$\chi_{\rm CF}'(G) \leq \log_2\Delta + C$.
\end{conjecture}

\begin{conjecture}\label{ConjectureWeak-closed}
For every graph $G$ of maximum degree $\Delta$ without isolated vertices,
$\chi_{\rm CF}'(G) \leq (1+o(1))\log_2\Delta$.
\end{conjecture}

\begin{conjecture}\label{ConjectureStrong-open}
There is a constant $C$ such that for every graph $G$ of maximum degree $\Delta$ without isolated vertices and edges,
$\chi_{\rm OCF}'(G) \leq \log_2\Delta + C$.
\end{conjecture}

\begin{conjecture}\label{ConjectureWeak-open}
For every graph $G$ of maximum degree $\Delta$ without isolated vertices and edges,
$\chi_{\rm OCF}'(G) \leq (1+o(1))\log_2\Delta$.
\end{conjecture}

\end{document}